\documentclass[preprint,12pt,authoryear]{elsarticle}

\usepackage{amssymb}

\usepackage{amsthm}

\usepackage{mathtools}
\newtheorem{thm}{Theorem}
\usepackage{multirow}
\usepackage{amsmath,amsfonts,amssymb,amsthm,epsfig,epstopdf,url,array}
\usepackage{caption}
\usepackage{subcaption}

\journal{Nuclear Physics B}

\begin{document}
	
	\begin{frontmatter}
		
		
		
		\tnotetext[label1]{corresponding author.}
		\author{Amenah AL-Najafi\corref{cor1}\fnref{label1}}
		\ead{amenah@math.u-szeged.hu}
		\author{L\'aszl\'o Viharos\corref{cor1}\fnref{label1}}
		\ead{viharos@math.u-szeged.hu}
		
		
		\title{Weighted least squares estimators for the Parzen tail index}
		
		\address[label1]{Bolyai Institute, University of Szeged, Aradi v\'ertan\'uk tere 1, Szeged, H-6720, Hungary}



\begin{abstract}
Estimation of the tail index of heavy-tailed distributions and its applications are essential in many research areas. We propose a class of weighted least squares (WLS) estimators for the Parzen tail index. Our approach is based on the method developed by \cite{Holan2010}. We investigate consistency and asymptotic normality of the WLS estimators. Through a simulation study, we make a comparison with the Hill, Pickands, DEdH (Dekkers, Einmahl and de Haan) and
ordinary least squares (OLS) estimators using the mean square error as criterion. The results show that in a restricted model some members of the WLS estimators are competitive with the Pickands, DEdH and OLS estimators.

\end{abstract}


\end{frontmatter}

\section{Introduction and main results}
The problem of estimating the tail characteristics of probability distributions has received enormous attention in the last decades.		
Let $F$ be an absolutely continuous probability distribution function with density function $f$ and let $Q$ denote the corresponding quantile function defined as
$$Q(s) := \inf\{x:F(x)\ge s\}, \;\; 0<s\le 1, \;\;\; Q(0) := Q(0+).$$
\cite{Parzen1979} used the density-quantile function $fQ(\cdot)=f(Q(\cdot))$ to classify probability distributions. Parzen assumed that the limit
\begin{equation}\nu_1:=\lim_{u\to1}\frac{(1-u)J(u)}{fQ(u)}\label{nu1}\end{equation}
exists, where $J$ is the score function defined as $J(u)=-(fQ)'(u)$.
Assumption (\ref{nu1}) yields the following approximation for $u$ values near 1:
$$fQ(u)\approx C(1-u)^{\nu_1},$$
for some positive constant $C$.
Based on the parameter $\nu_1$, Parzen classified the probability distributions. Heavy tailed distributions correspond to $\nu_1>1$.

\cite{Parzen2004} assumed that $fQ(\cdot)$ is regularly varying at 0 and 1:
\begin{align}
fQ(u)&=u^{\nu_0}L_0(u),\quad u\in[0,1/2),\label{lefttail}\\
fQ(u)&=(1-u)^{\nu_1}L_1(1-u),\quad u\in(1/2,1],\label{righttail}
\end{align}
where $\nu_0,\nu_1>0$ are finite constants and $L_0$ and $L_1$ are slowly varying at zero.
The parameters $\nu_0$ and $\nu_1$ are called the left and right tail exponents of the density-quantile function.

\cite{Holan2010} considered the following orthogonal series expansion for $L_i$:
\begin{equation}L_i(u)=\exp \bigg\lbrace \theta_{i,0}+2\sum_{k=1}^{\infty}\theta_{i,k}\cos(2\pi ku)\bigg\rbrace,\quad i=0,1.\label{orthogonal}\end{equation}
In order to estimate the tail exponents, they assumed that $L_i$ admits the representation
\[L_i(u)=L_i^{(p_i)}(u)=\exp \bigg\lbrace \theta_{i,0}+2\sum_{k=1}^{p_i}\theta_{i,k}\cos(2\pi ku)\bigg\rbrace,\quad i=0,1,\]
where $p_i$ is fixed and unknown. In the representation (\ref{lefttail}) and (\ref{righttail}) they considered $fQ(u)$ for $u\in(0,u_l]$ and $u\in[u_r,1)$, where $u_l\le1/2$ and $u_r\ge1/2$ are chosen by the statistician. Using these representations, they assumed the equations
\begin{align*}
\log fQ(u)&=\nu_0\log u+\theta_{0,0}+2\sum_{k=1}^{p_0}\theta_{0,k}\cos(2\pi ku),\quad u\in(0,u_l],\\
\log fQ(u)&=\nu_1\log(1-u)+\theta_{1,0}+2\sum_{k=1}^{p_1}\theta_{1,k}\cos(2\pi k(1-u)),\quad u\in[u_r,1).
\end{align*}
Based on some estimator $\widehat{fQ}(u)$ of the density-quantile $fQ(u)$, this leads to the regression equations
\begin{align*}
\log\widehat{fQ}(u_j)&=\nu_0\log u_j+\theta_{0,0}+2\sum_{k=1}^{p_0}\theta_{0,k}\cos(2\pi ku_j)+\varepsilon(u_j),\\
\log\widehat{fQ}(1-u_j)&=\nu_1\log u_j+\theta_{1,0}+2\sum_{k=1}^{p_1}\theta_{1,k}\cos(2\pi k u_j)+\varepsilon(1-u_j),
\end{align*}
where $\varepsilon(u)=\log\big(\widehat{fQ}(u)/fQ(u)\big)$ is the residual process, $u_{j}=j/n, j=u_{\lceil na\rceil},\dots, u_{\lfloor nb\rfloor}$ and $0<a<b<1$, so the percentiles $u_{j}$ are chosen from a subset $U=[a,b]$ of the interval $(0,1).$
 \cite{Holan2010} obtained some estimators $\widehat\nu_0$ and $\widehat\nu_1$ for the tail exponents $\nu_0$ and $\nu_1$ using ordinary least squares regression. 

We propose a more general class of estimators using weighted least squares regression. We choose some nonnegative weights of the form $w_{j,n}=R(j/n)$ with some weight function $R$. Set $y_j:=\log\widehat{fQ}(u_j)$,
\begin{align*}
y&:=(y_{\lceil na\rceil},\ldots,y_{\lfloor nb\rfloor})',\\
W&:=\textrm{diag}(w_{\lceil na\rceil,n},\ldots,w_{\lfloor nb\rfloor,n}),
\end{align*}
and let $X:=[G^*,G_0,2G_1,\ldots,2G_{\widetilde p_0}]$, where
\begin{align*}
G^*&=\big(\log(u_{\lceil na\rceil}),\ldots,\log(u_{\lfloor nb\rfloor})\big)'\\
G_{k}&=\big(\cos(2\pi ku_{\lceil na\rceil}),\ldots,\cos(2\pi ku_{\lfloor nb\rfloor})\big)',\quad k=0,\ldots,\widetilde p_0,
\end{align*}
and $\widetilde p_0>p_0$ is chosen by the statistician. Set $\beta_{\widetilde p_0}:=(\nu_0,\theta_{0,0},\theta_{0,1},\ldots,\theta_{0,\widetilde p_0})'$, where $\theta_{0,j}=0$ if $j>p_0$. By minimizing the weighted sum of squares
\[\sum_{j=\lceil na\rceil}^{\lfloor nb\rfloor}w_{j,n}\big(y_j-\nu_0\log u_j-\theta_{0,0}-2\sum_{k=1}^{\widetilde p_0}\theta_{0,k}\cos(2\pi ku_j)\big)^2,\]
we obtain the following estimator of $\beta_{\widetilde p_0}$:
\[\widehat{\beta}_{\widetilde p_0}=(X'WX)^{-1}X'Wy.\]
Then the weighted least squares estimator of $\nu_0$ can be written in the form
\[\widehat \nu_0=e_1'\widehat{\beta}_{\widetilde p_0}=e_1'(X'WX)^{-1}X'Wy,\]
where $e_1$ is the $\widetilde p_0+2$ dimensional vector defined as $e_1=(1,0,0,\ldots,0)'$. The right tail exponent $\nu_1$ can be estimated similarly.

A crucial point of this method is to choose a good estimator for the density-quantile $fQ(u)$. Letting $q(u):=Q'(u)$ denote the quantile density function, and using the identity $fQ(u)Q'(u)=1$, one wish to estimate $q(u)$ instead of $fQ(u)$. Given a sample $X_1, \ldots, X_n$ with distribution function $F$, let $F_n$ denote its empirical distribution function and define $Q_n:=F_n^{-1}$ to be the empirical quantile function. \cite{Holan2010} used the kernel quantile estimator of $q(u)$:
\begin{equation}\label{kernelest}
\widehat q_n(u)=\frac d{du}\int_0^1Q_n(t)K_n(u,t)d\mu_n(t),\quad u\in(0,1),
\end{equation} 
where the kernel function $K_n(u,t)$ and the measure $\mu_n$ satisfy the following conditions of \cite{Cheng1995}:

\bigskip\noindent
$(K_1)$ For every $n$, $0<\mu_n ([0,1])<\infty$, and $\mu_n(\left\lbrace 0,1\right\rbrace )=0$.\\
$(K_2)$ For every $n$ and each $(u,t)$, $K_n(u,t) \ge0$, and for every $u\in U$, $\int_0^1K_n(u,t)d \mu_n(t)=1$.\\
$(K_3)$ For every $n$, $\int_0^1tK_n(u,t)d\mu_n(t)=u, u\in U$.\\
$(K_4)$ There is a sequence $\delta_n\downarrow 0$ such that $\sup_{u\in U}\big|\int_{u-\delta_n}^{u+\delta_n}K_n(u,t)d\mu_n(t)-1\big|\downarrow0$ as $n\uparrow\infty$.

\bigskip
Let $S_n$ be the unique closed subset of (0,1) such that $\mu_n\big((0,1)\backslash S_n\big) =0$ and $\mu_n\big((0,1)\backslash S'_n\big)>0$ for any $S'_n\subset S_n$. For the sequence $\delta_n$ in $(K_4)$, let $I_n(u)=[u-\delta_n,u+\delta_n]$, $I_n^{c}(u)=(0,1)\backslash I_n(u),$ for $u\in U$. Define $\Lambda(u;K_n)=\int_{I_n(u)}|K'_n(u,t)|d\mu_n(t),$ $u\in U$, and for a well-defined  function $g$ on (0,1), let $\Psi(g;K_n)=\sup_{u\in U}\int_{I_n^c(u)}|g(t)K'_n(u,t)|d \mu_n(t)$. It is also assumed that the derivative $K'_n(u,t)=\partial K_n(u,t)/\partial u$ satisfies the conditions $(K_5)-(K_7)$ below: 

\bigskip\noindent
$(K_5)$ For every $n$, $\sup_{u\in U}\int_0^1|K'_n(u,t)|d\mu_n(t)<\infty$.\\
$(K_6)$ For every $n$ and each $u\in U$, $K_n(u,t)\equiv 0$, $t\in I_n^{c}(u)$; or $S_n\subseteq [\varepsilon,1-\varepsilon]\subset(0,1),$ with $U\subset [\varepsilon,1-\varepsilon]$ for some $0<\varepsilon<1/2$.\\
$(K_7)$ For the sequence $\delta_n$ in $(K_4)$, $\delta_n^2\sup_{u\in U}\Lambda (u;K_n)\rightarrow 0$ and $\Psi(1;K_n)\to0$ as $n\uparrow\infty$.

\bigskip
Similarly as in \cite{Holan2010}, in some cases we assume that the kernel function has the form $K_n(u,t)=K\big(h_n^{-1}(t-u)\big)h_n^{-1}$ and satisfies the condition
\[(K_8)\quad\sup_{u\in U}\left|h_n^{-1}K\Big(\frac{s-u}{h_n}\Big)-h_n^{-1}K\Big(\frac{t-u}{h_n}\Big)\right|\le C_n|t-s|^\beta\quad {\rm and}\quad |K''(x)|\le C/|x|\] 
for some constants $C,\beta>0$ and $|x|$ sufficiently large, and $C_n$ are positive constants such that $\sup_{n\ge1}C_n<\infty$.

\bigskip
Moreover, \cite{Holan2010} used the following assumptions of \cite{Cheng1995} on $q(u)$:

\bigskip\noindent
$(Q_1)$ The quantile density function is twice differentiable on (0,1).\\
$(Q_2)$ There exists a positive constant $\gamma$ such that $\sup_{u\in (0,1)}u(1-u)|J(u)|/fQ(u)\\ \leq \gamma$, where $J$ is the score function in (\ref{nu1}).\\
$(Q_3)$ Either $q(0)<\infty $ or $q(u)$ is nonincreasing in some interval $(0,u_{*})$, and either $q(1)<\infty $ or $q(u)$ is nondecreasing in some interval $(u^*,1)$.\\

We will show that the limit matrix $ M(a,b,R):=\lim_{n\to\infty}n^{-1}X'WX$ exists (see the proof of Theorem \ref{t1}). Let $(v^*, v_0,\ldots, v_{\widetilde p_i})$ be the first row of $M(a,b,R)^{-1}$, and set $G_R(u):=R(u)\big(v^*\log u+v_0+2\sum_{k=1}^{\widetilde p_i}v_k\cos(2\pi ku)\big)$, $i=0,1$.

\bigskip
Finally, we assume that the weight function $R$ satisfies the following condition:

\bigskip\noindent
$(R)$ $R$ is nonnegative and Riemann integrable on $[a,b]$.

\bigskip
Let $\overset P\longrightarrow$ denote convergence in probability, $\overset{\cal D}\longrightarrow$ denote convergence in distribution, and let $N(\mu,\sigma^2)$
stand for the normal distribution with mean $\mu$ and variance
$\sigma^2$.
Limiting and order relations are always meant as $n\to\infty$
if not specified otherwise.
Our main results are contained in the following two theorems:

\begin{thm}\label{t1}
	Suppose that the conditions $(Q_1)-(Q_3)$ are satisfied for the quantile density $q(u)$, and $\widehat q(u)$ is a kernel smoothed estimator with kernel function satisfying $(K_1)-(K_7)$, the weight function $R$ satisfies the condition $(R)$, and the matrix $M(a,b,R)$ is invertible. Moreover, assume that the percentiles $u_j$ are chosen from a closed set $U=[a,b]$ such that $u_j=j/n$, $j=\lceil na\rceil,\ldots,\lfloor nb\rfloor$, and $\tilde p_i>p_i$, $i=0,1$. Then $\widehat \nu_i\overset P\longrightarrow\nu_i$, $i=0,1$.
\end{thm}

\begin{thm}\label{t2}
	Assume that the conditions of Theorem \ref{t1} are satisfied, and suppose that the kernel function is symmetric and differentiable on $[-1,1]$, and satisfies the condition $(K_8)$. Suppose that the derivative $g_R(u):=G_R'(u)$ exists, and $g_R$ and $G_R$ are uniformly bounded on $U$. Let $h_n$ be a sequence such that $nh^2_n\rightarrow\infty$, $nh^4_n\rightarrow0$ and $h_n\rightarrow0$, and assume that $\tilde p_i>p_i$, $i=0,1$. Then
	\[\sqrt{n}(\widehat\nu_i-\nu_i)\overset{\cal D}\longrightarrow N(0,V),\quad i=0,1,\]
	where
	\begin{equation}
	\label{variance}
	V=\int_a^bG^2_R(u)du+\int_a^b\int_a^bG_R(u)G_R(v)\bigg(1+[(u\land v)-uv]\dfrac{q'(u)q'(v)}{q(u)q(v)} \bigg)dudv.
	\end{equation}
\end{thm}

\bigskip
In the special case when the weight function $R$ is identically 1, the two theorems above reduces to Theorems 1 and 2 of \cite{Holan2010}.

\section{Classical tail index estimation}

A distribution function $F$ has right heavy tail with tail index $\alpha_1>0$ if $1-F(x)$ is regularly varying at infinity with index $-1/\alpha_1$, i.e.,
\begin{equation}
\label{Pareto_tail}
1-F(x)=x^{-1/\alpha_1}\ell_1(x),\quad 0<x<\infty,
\end{equation}
where $\ell_1$ is a function slowly varying at infinity. Similarly, $F$ has left heavy tail with tail index $\alpha_0>0$ if $F(-x)$ is regularly varying at infinity with index $-1/\alpha_0$.
Let $X_1, X_2,\ldots$ be independent random variables with a common distribution function $F$ having right heavy tail with tail index $\alpha_1$, and for each $n \in \mathbb N$, let $X_{1,n}\le \cdots \le X_{n,n}$ denote the order statistics pertaining to the sample $X_1, \ldots, X_n$.
Several estimators exist for $\alpha_1$ among which Hill's estimator
is the most classical. \cite{Hill1975} proposed the following estimator for the tail index $\alpha_1$:
\[ \widehat\alpha_1=\frac1{k_n}\sum_{j=1}^{k_n}\log X_{n-j+1,n}-\log X_{n-k_n,n}=\frac1{k_n}\sum_{j=1}^{k_n}\log\frac{X_{n-j+1,n}}{X_{n-k_n,n}},\]
where the $k_n$ are some integers satisfying
$$1\le k_n < n, \quad\ k_n\to\infty\quad\ {\rm and}\quad\ k_n/n\to0\ \quad{\rm as}\ n\to\infty.$$
The left tail analogue of the Hill estimator is the following:
\[ \widehat\alpha_0=\frac1{k_n}\sum_{j=1}^{k_n}\log\frac{X_{j,n}}{X_{k_n+1,n}}.\]
Another estimators were proposed by \cite{Pickands1975}, \cite{Dekkers1989}, to name a few.
As it was pointed out by \cite{Holan2010}, there is the following relationship between the Parzen and classical tail indices: $\nu_i=1+\alpha_i$, $i=0,1$. Thus the classical tail index estimators also can be used to estimate the Parzen tail index.

\section{Comparison of tail index estimators}

\subsection{Asymptotic variances}

We evaluate the limiting variance (\ref{variance}) for $\widetilde p_0=1$, different weight functions and tail indices to compare the WLS and the unweighted (ordinary least squares) estimators in the following submodel of (\ref{orthogonal}):
\[L_0(u)=\exp\big\lbrace 2\cos(2\pi u)\big\rbrace,\quad u\in[a,b].\]
The limiting variances are contained in Table \ref{table_variance} in the Appendix. For the calculations we used numerical integration performed by the Wolfram Mathematica software. We see that in some cases the use of the weights makes the asymptotic variance smaller.

\subsection{Simulation results}

In order to make a comparison with existing proposals, simulations were done performed by the Matlab software. The samples were generated from the model (\ref{lefttail}) with $L_0\equiv1$ using different tail indices $\nu_0$. The Hill, Pickands, DEdH (Dekkers, Einmahl and de Haan) and the least squares estimators were included in the simulation study. Similarly as in \cite{Holan2010}, for the simulations we used the Bernstein polynomial estimator of $q(u)$. Let $0<\varepsilon<1/2$ be a constant, and assume that $U\subset[\varepsilon,1-\varepsilon]$. Set $L_\varepsilon:=1-2\varepsilon$ and $t_j:=\varepsilon+(j/k)L_\varepsilon$, $j=0,1,\ldots,k$. The Bernstein polynomial estimator is defined as
\[\widehat{q}^{B}_n(u)=\dfrac 1{L^{k}_\varepsilon}\sum_{j=0}^{k-1}\frac{Q_n(t_{j+1})-Q_n(t_j)}{1/k}\binom{k-1}j(u-\varepsilon)^j(1-\varepsilon-u)^{k-1-j}.\]
This estimator belongs to the class (\ref{kernelest}) and satisfies the conditions $(K_1)-(K_7)$.
We used the values $k=n=700$, $\varepsilon=0.001$, $a=0.001$ and $b=0.4$ for the regression estimators, and the weight function $R(u)=u/300$ for the WLS estimator.
Tables \ref{table_mean} and \ref{table_mse} contain the average simulated estimates (mean) and the calculated empirical mean square errors (MSE). We used the sample fraction size $k_n=100$ for the Hill, Pickands and DEdH estimators. All the simulations were repeated 200 times.
We conclude that in the submodel $L_0\equiv1$ for $\alpha$ values between 0.8 and 1.5 the WLS estimator has better performance than the OLS estimator. Thus for thinner tails we propose the WLS estimator instead of the OLS estimator. The Hill estimator is the best among the examined estimators. This good performance is not surprising since the Hill estimator was obtained in the special case of (\ref{Pareto_tail}) when the slowly varying function $\ell_1(x)$ is constant for all $x\ge x_{\alpha_1}$, for some threshold $x_{\alpha_1}$. The Pickands estimator has also good performance. On the other hand, we emphasize that the WLS method can be applied not only for the estimation of the tail index but for the estimation of the slowly varying functions $L_i$ in (2) and (3).

\section{Proofs}

The proof of Theorems \ref{t1} and \ref{t2} follows the general outline of the proof of Theorems 1 and 2 of \cite{Holan2010}. We give a more detailed proof for Theorem \ref{t2}.

\bigskip\noindent
{\bf Proof of Theorem \ref{t1}.} We deal only with the left tail exponent $\nu_0$, the proof for $\nu_1$ is similar. Set $\gamma=(\gamma_{\lceil na\rceil},\ldots,\gamma_{\lfloor nb\rfloor})':=\sqrt{W}X(X'WX)^{-1}e_{1}$ and $\underline\varepsilon:=\big(\varepsilon(u_{\lceil na\rceil}),\ldots,\varepsilon(u_{\lfloor nb\rfloor})\big)'$. Then $\widehat\nu_0-\nu_0=\gamma'\sqrt W\underline\varepsilon$, and hence, using the Cauchy-Schwarz inequality,
\begin{align*}
|\widehat\nu_0-\nu_0|=\Bigg|\sum_{j=\lceil na\rceil}^{\lfloor nb\rfloor}\gamma_j\sqrt{w_{j,n}}\,\varepsilon(u_j)\Bigg|
\le\Bigg(\sum_{j=\lceil na\rceil}^{\lfloor nb\rfloor}\gamma_j^2\Bigg)^{1/2}\Bigg(\sum_{j=\lceil na\rceil}^{\lfloor nb\rfloor}w_{j,n}\,\varepsilon^2(u_j)\Bigg)^{1/2}.
\end{align*}
We have $\sum_{j=\lceil na\rceil}^{\lfloor nb\rfloor}\gamma^2_j=\gamma'\gamma=e'_1(n^{-1}X'WX)^{-1}e_1n^{-1}$ with the matrix
\begin{small}
	\[X'WX=\begin{bmatrix}
	\smashoperator\sum_{j=\lceil na\rceil}^{\lfloor nb\rfloor}\log^2u_jR(u_j)
	& \smashoperator\sum_{j=\lceil na\rceil}^{\lfloor nb\rfloor}\log u_jR(u_j) & 2\smashoperator\sum_{j=\lceil na\rceil}^{\lfloor nb\rfloor}\log u_j\cos(2\pi u_j)R(u_j)\ldots\\
	\smashoperator \sum_{j=\lceil na\rceil}^{\lfloor nb\rfloor}\log u_jR(u_j) & \sum_{j=\lceil na\rceil}^{\lfloor nb\rfloor}R(u_j) & 2\smashoperator\sum_{ j=\lceil na\rceil}^{\lfloor nb\rfloor}\cos(2\pi u_j)R(u_j)\ldots\\
	\vdots & \vdots & \vdots \end{bmatrix}.\]
\end{small}	
Then by Riemann sum approximation
\begin{align}\label{mlimit}
\lim_{n\to\infty}&n^{-1}X'WX=M(a,b,R)\\&:=
\begin{bmatrix}
\int_a^b\log^2u\,R(u)du
& \int_a^b\log u\,R(u)du & 2\int_a^b\log u\cos(2\pi u)R(u)du\ldots\\
\int_a^b\log u\,R(u)du & \int_a^bR(u)du & 2\int_a^b\cos(2\pi u)R(u)du\ldots\\
\vdots & \vdots & \vdots \end{bmatrix}.\nonumber
\end{align}
It follows that for all $n$ large enough $e'_1(n^{-1}X'WX)^{-1}e_1\le C$ for some constant $C$, and hence
\[|\widehat\nu_0-\nu_0|\le\sqrt C\bigg(n^{-1}\sum_{j=\lceil na\rceil}^{\lfloor nb\rfloor}R(u_j)\,\varepsilon^2(u_j)\bigg)^{1/2}.\]
Let $C'>0$ be a constant such that $R(u)\le C'$, $0\le u\le1$. Then
\[|\widehat\nu_0-\nu_0|\le\sqrt{CC'}\bigg(n^{-1}\sum_{j=\lceil na\rceil}^{\lfloor nb\rfloor}\varepsilon^2(u_j)\bigg)^{1/2}.\]
Now, by Theorem 2.1 of \cite{Cheng1995}, $n^{-1}\sum_{j=\lceil na\rceil}^{\lfloor nb\rfloor}\varepsilon^2(u_j)=o_P(1)$ (cf.~the proof of Theorem 1 of \cite{Holan2010}). \quad$\square$

\bigskip\noindent
{\bf Proof of Theorem \ref{t2}.}
Write
\begin{align*}
\sqrt n(\widehat{v_0}-v_0)&=\frac1{\sqrt{n}}e'_1(n^{-1}X'WX)^{-1}X'W\underline\varepsilon\\
&=\frac1{\sqrt n}e'_1M(a,b,R)^{-1}X'W\underline\varepsilon+\frac1{\sqrt n}e'_1\big((n^{-1}X'WX)^{-1}-M(a,b,R)^{-1}\big)X'W\underline\varepsilon.
\end{align*}
By straightforward calculation,
\begin{equation}\label{eproc}
A_n:=\frac1{\sqrt n}e'_1M(a,b,R)^{-1}X'W\underline\varepsilon=\frac1{\sqrt n}\sum_{j=\lceil na\rceil}^{\lfloor nb\rfloor}\varepsilon(u_j)G_R(u_j).\end{equation}
It follows from Theorem 5 in \cite{Holan2010} that
\[A_n\overset{\cal D}\longrightarrow G_R(b)W(b)-G_R(a)W(a)-\int_a^bW(u)\left(g_R(u)-G_R(u)\dfrac{q'(u)}{q(u)}\right)du,\]
where $W(u)$ is a Brownian bridge process. The limiting variance is given by (\ref{variance}). Next we show that
\[B_n:=\frac1{\sqrt n}e'_1\big((n^{-1}X'WX)^{-1}-M(a,b,R)^{-1}\big)X'W\underline\varepsilon=o_P(1).\]
Let $(v^*_n, v_{0,n},\ldots, v_{\widetilde p_0,n})$ be the first row of $(n^{-1}X'WX)^{-1}-M(a,b,R)^{-1}$.
By (\ref{mlimit}),
$(v^*_n, v_{0,n},\ldots, v_{\widetilde p_0,n})\to\,\mathbf0$.
Set
\[G^{(n)}(u)=R(u)\big(v^*_n\log u+v_{0,n}+2\sum_{k=1}^{\widetilde p_0}v_{k,n}\cos(2\pi ku)\big).\]
Similarly as in (\ref{eproc}),
\begin{align*}
B_n&=\dfrac1{\sqrt{n}}\sum_{j=\lceil na\rceil}^{ \lfloor nb\rfloor}\varepsilon(u_j)G^{(n)}(u_j)\\
&=v^*_n\dfrac1{\sqrt{n}}\sum_{j=\lceil na\rceil}^{ \lfloor nb\rfloor}\varepsilon(u_j)R(u_j)\log u_j+v_{0,n}\dfrac1{\sqrt{n}}\sum_{j=\lceil na\rceil}^{ \lfloor nb\rfloor}\varepsilon(u_j)R(u_j)\\
&\ \ +2\sum_{k=1}^{\widetilde p_0}v_{k,n}\dfrac1{\sqrt{n}}\sum_{j=\lceil na\rceil}^{ \lfloor nb\rfloor}\varepsilon(u_j)R(u_j)\cos(2\pi ku_j).
\end{align*}
Each term in the last sum tends to zero, e.g., in the first term $v^*_n\to0$ and using again Theorem 5 in \cite{Holan2010}, the sequence $\frac1{\sqrt n}\sum_{j=\lceil na\rceil}^{ \lfloor nb\rfloor}\varepsilon(u_j)R(u_j)\log u_j$ has a weak limit.
\quad$\square$

\vfill\eject
\section*{Appendix}

\begin{table}[!h]
	\centering
	\caption{Limiting variances for different weight functions and tail indices.}
	\label{table_variance}
	\begin{tabular}{l|llll|l}
		\hline
		\multirow{3}{*}{$\nu_0=1.2$} & \multicolumn{4}{c|}{\multirow{2}{*}{R(u)}}  & \multicolumn{1}{l}{\multirow{3}{*}{unweighted}} \\
		& \multicolumn{4}{c|}{}                          & \multicolumn{1}{c}{}                            \\ \cline{2-5}
		& $1+\cos u$         & $e^{-u}$ & $-\log u$        &$1/u$& \multicolumn{1}{c}{}                           \\ \hline
		$a=0.1, b=0.4$                  & 821.232         & 816.812   & 823.778 &851.364& 822.13                                           \\
		$a=0.1, b=0.3$                  & 1512.62         & 1513.46    &1538.35 & 1600.46& 1512.83                                         \\
		$a=0.2, b=0.3$                  & 269523 & 269655    & 270796 &272081       & 269524                                          \\ \hline
	\end{tabular}
\end{table}

\begin{table}[!h]
	\centering
	\begin{tabular}{l|llll|l}
		\hline
		\multirow{3}{*}{$\nu_0=1.8$} & \multicolumn{4}{c|}{\multirow{2}{*}{R(u)}}  & \multicolumn{1}{l}{\multirow{3}{*}{unweighted}} \\
		& \multicolumn{4}{c|}{}                          & \multicolumn{1}{c}{}                            \\ \cline{2-5}
		& $1+\cos u$         & $e^{-u}$ & $-\log u$        &$1/u$& \multicolumn{1}{c}{}                            \\ \hline
		$a=0.1, b=0.4$                  & 821.962         & 819.166   & 829.786 &860.498& 822.66                                           \\
		$a=0.1, b=0.3$                  & 1521.58         & 1523.69    &1551.68 & 1617.04& 1521.66                                        \\
		$a=0.2, b=0.3$                  & 267666 & 267807    & 268969 &270267      & 267666                                          \\ \hline
	\end{tabular}
\end{table}

\begin{table}[!h]
	\centering
	\begin{tabular}{l|llll|l}
		\hline
		\multirow{3}{*}{$\nu_0=1.667$} & \multicolumn{4}{c|}{\multirow{2}{*}{R(u)}}  & \multicolumn{1}{c}{\multirow{3}{*}{unweighted}} \\
		& \multicolumn{4}{c|}{}                          & \multicolumn{1}{c}{}                            \\ \cline{2-5}
		& $1+\cos u$         & $e^{-u}$ & $-\log u$        &$1/u$& \multicolumn{1}{c}{}                            \\ \hline
		$a=0.1, b=0.4$                  & 819.423         & 816.278   & 826.109 &856.14& 820.164                                          \\
		$a=0.1, b=0.3$                  & 1516.49         & 1518.31    &1545.6 & 1610.22& 1516.6                                        \\
		$a=0.2, b=0.3$                  & 268011 & 268151    & 269308 &270604     & 268012                                         \\ \hline
	\end{tabular}
\end{table}

\begin{table}[!h]
	\centering
	\begin{tabular}{l|llll|l}
		\hline
		\multirow{3}{*}{$\nu_0=2.25$} & \multicolumn{4}{c|}{\multirow{2}{*}{R(u)}}  & \multicolumn{1}{c}{\multirow{3}{*}{unweighted}} \\
		& \multicolumn{4}{c|}{}                          & \multicolumn{1}{c}{}                            \\ \cline{2-5}
		& $1+\cos u$         & $e^{-u}$ & $-\log u$        &$1/u$& \multicolumn{1}{c}{}                            \\ \hline
		$a=0.1, b=0.4$                  & 840.595         & 838.929   & 825.157 &885.102& 841.151                                           \\
		$a=0.1, b=0.3$                  & 1551.91        & 1555.02    &1585.51 & 1653.45& 1551.89                                      \\
		$a=0.2, b=0.3$                  & 266776 & 266924    & 268099 &269406      & 266775                                         \\ \hline
	\end{tabular}
\end{table}

\begin{table}[!h]
	\centering
	\caption{Average simulated tail index estimates (Mean) for sample size $n=700$ and for $L_0\equiv1$.}
	\label{table_mean}
	\scalebox{0.85}{
		\begin{tabular}{cccccccccc}
			\hline
			\multirow{2}{*}{}          & \multicolumn{9}{c}{Mean}                                                                                                                                                                                            \\ \cline{2-10} 
			& \multicolumn{3}{c|}{WLS}                    & \multicolumn{3}{c|}{OLS}                     & \multicolumn{1}{c|}{\multirow{2}{*}{Hill}} & \multicolumn{1}{c|}{\multirow{2}{*}{Pickands}} & \multirow{2}{*}{DEdH} \\ \cline{1-7}
			\multicolumn{1}{c|}{$\nu(\alpha)$}     & $\widetilde p_0=1$    & $\widetilde p_0=2$      & \multicolumn{1}{c|}{$\widetilde p_0=3$}      & $\widetilde p_0=1$    & $\widetilde p_0=2$      & \multicolumn{1}{c|}{$\widetilde p_0=3$}      & \multicolumn{1}{c|}{}                      & \multicolumn{1}{c|}{}                          &                       \\ \hline
			\multicolumn{1}{c|}{2.25(1.25)}  & 2.3777 & 2.4751 & \multicolumn{1}{c|}{2.5088} & 2.4271 & 2.4803 & \multicolumn{1}{c|}{2.4825} & \multicolumn{1}{c|}{2.2396}                & \multicolumn{1}{c|}{2.2703}                    & 2.7346                \\
			\multicolumn{1}{c|}{2(1)}     & 2.0741 & 2.1231 & \multicolumn{1}{c|}{2.2423} & 2.0902 & 2.1162 & \multicolumn{1}{c|}{2.1177} & \multicolumn{1}{c|}{2.0038}                & \multicolumn{1}{c|}{1.9998}                    & 2.4988                \\
			\multicolumn{1}{c|}{1.833(0.833)} & 1.9119 & 1.9249 & \multicolumn{1}{c|}{1.9405} & 1.9248 & 1.904  & \multicolumn{1}{c|}{1.8959} & \multicolumn{1}{c|}{1.8404}                & \multicolumn{1}{c|}{1.8471}                    & 2.3354                \\
			\multicolumn{1}{c|}{1.667(0.667)} & 1.7163 & 1.6915 & \multicolumn{1}{c|}{1.7274} & 1.7217 & 1.7019 & \multicolumn{1}{c|}{1.7058} & \multicolumn{1}{c|}{1.6743}                & \multicolumn{1}{c|}{1.6902}                    & 2.1692                \\
			\multicolumn{1}{c|}{1.556(0.556)} & 1.5949 & 1.6294 & \multicolumn{1}{c|}{1.5951} & 1.6017 & 1.5822 & \multicolumn{1}{c|}{1.5637} & \multicolumn{1}{c|}{1.5534}                & \multicolumn{1}{c|}{1.5567}                    & 2.0483                \\
			\multicolumn{1}{c|}{1.5(0.5)}   & 1.5239 & 1.5448 & \multicolumn{1}{c|}{1.5518} & 1.5222 & 1.5613 & \multicolumn{1}{c|}{1.5668} & \multicolumn{1}{c|}{1.5005}                & \multicolumn{1}{c|}{1.4942}                    & 1.9955                \\
			\multicolumn{1}{c|}{1.333(0.333)} & 1.3639 & 1.389  & \multicolumn{1}{c|}{1.3874} & 1.3598 & 1.3335 & \multicolumn{1}{c|}{1.3136} & \multicolumn{1}{c|}{1.3347}                & \multicolumn{1}{c|}{1.3294}                    & 1.8296                \\
			\multicolumn{1}{c|}{1.25(0.25)}  & 1.2956 & 1.2471 & \multicolumn{1}{c|}{1.242}  & 1.2741 & 1.2585 & \multicolumn{1}{c|}{1.2629} & \multicolumn{1}{c|}{1.2476}                & \multicolumn{1}{c|}{1.2474}                    & 1.7426                \\
			\multicolumn{1}{c|}{1.2(0.2)}   & 1.2281 & 1.2483 & \multicolumn{1}{c|}{1.2189} & 1.1967 & 1.2204 & \multicolumn{1}{c|}{1.2089} & \multicolumn{1}{c|}{1.1993}                & \multicolumn{1}{c|}{1.2144}                    & 1.6942                \\
			\multicolumn{1}{c|}{1.182(0.182)} & 1.1742 & 1.1891 & \multicolumn{1}{c|}{1.199}  & 1.1776 & 1.1725 & \multicolumn{1}{c|}{1.1677} & \multicolumn{1}{c|}{1.1833}                & \multicolumn{1}{c|}{1.174}                     & 1.6783                \\
			\multicolumn{1}{c|}{1.167(0.167)} & 1.1628 & 1.1953 & \multicolumn{1}{c|}{1.1826} & 1.162  & 1.158  & \multicolumn{1}{c|}{1.1452} & \multicolumn{1}{c|}{1.167}                 & \multicolumn{1}{c|}{1.1624}                    & 1.662                 \\
			\multicolumn{1}{c|}{1.1(0.1)}   & 1.1116 & 1.0926 & \multicolumn{1}{c|}{1.1538} & 1.0899 & 1.0755 & \multicolumn{1}{c|}{1.0725} & \multicolumn{1}{c|}{1.1006}                & \multicolumn{1}{c|}{1.0952}                    & 1.5955                \\
			\multicolumn{1}{c|}{1.067(0.067)} & 1.0761 & 1.106  & \multicolumn{1}{c|}{1.0895} & 1.0456 & 1.0597 & \multicolumn{1}{c|}{1.0431} & \multicolumn{1}{c|}{1.0673}                & \multicolumn{1}{c|}{1.0562}                    & 1.5622                \\
			\multicolumn{1}{c|}{1.05(0.05)}  & 1.0674 & 1.0607 & \multicolumn{1}{c|}{1.0866} & 1.0527 & 1.0476 & \multicolumn{1}{c|}{1.0438} & \multicolumn{1}{c|}{1.0496}                & \multicolumn{1}{c|}{1.048}                     & 1.5445     \\ \hline             
	\end{tabular}}
\end{table}

\begin{table}[!h]
	\centering
	\caption{Empirical mean square errors (MSE) of tail index estimates for sample size $n=700$ and for $L_0\equiv1$.}
	\label{table_mse}
	
	\scalebox{0.82}{
		\begin{tabular}{ccccccclll}
			\hline
			\multicolumn{10}{c}{MSE}                                                                                                                                                                                                                                                                                                                  \\ \cline{2-10} 
			& \multicolumn{3}{c|}{WLS}                                                  & \multicolumn{3}{c|}{OLS}                                       & \multicolumn{1}{c|}{\multirow{2}{*}{Hill}} & \multicolumn{1}{c|}{\multirow{2}{*}{Pickands}} & \multicolumn{1}{c}{\multirow{2}{*}{DEdH}} \\ \cline{1-7}
			\multicolumn{1}{c|}{$\nu(\alpha)$}   & $\widetilde p_0=1$                         & $\widetilde p_0=2$                           & \multicolumn{1}{c|}{$\widetilde p_0=3$}        & $\widetilde p_0=1$    & $\widetilde p_0=2$                          & \multicolumn{1}{c|}{$\widetilde p_0=3$}      & \multicolumn{1}{c|}{}                      & \multicolumn{1}{c|}{}                          & \multicolumn{1}{c}{}                      \\ \hline
			\multicolumn{1}{c|}{2.25(1.25)}   & 0.0953                     & 0.1565                    & \multicolumn{1}{c|}{0.2224}  & 0.1540 & 0.2701                     & \multicolumn{1}{c|}{0.3855} & \multicolumn{1}{l|}{0.0177874}             & \multicolumn{1}{l|}{0.0592}                    & 0.2525                                    \\
			
			\multicolumn{1}{c|}{2(1)}       & 0.0794                     & 0.1121                     & \multicolumn{1}{c|}{0.1865}  & 0.1029 & 0.1244                     & \multicolumn{1}{c|}{0.1942} & \multicolumn{1}{l|}{0.0112351}             & \multicolumn{1}{l|}{0.0491}                    & 0.2600                                    \\
			
			\multicolumn{1}{c|}{1.833(0.833)} & 0.0599                     & 0.1134                     & \multicolumn{1}{c|}{0.1550}  & 0.0714 & 0.1257                     & \multicolumn{1}{c|}{0.1673} & \multicolumn{1}{l|}{0.0075016}             & \multicolumn{1}{l|}{0.0427}                    & 0.2598                                    \\
			
			\multicolumn{1}{c|}{1.667(0.667)} & 0.0594                     & 0.0817                     & \multicolumn{1}{c|}{0.1164}  & 0.0565 & 0.0832                     & \multicolumn{1}{c|}{0.1218} & \multicolumn{1}{l|}{0.0062222}             & \multicolumn{1}{l|}{0.0412}                    & 0.2471                                    \\
			
			\multicolumn{1}{c|}{1.556(0.556)} & 0.0515                     & 0.0935                    & \multicolumn{1}{c|}{0.0938}  & 0.0404 & 0.0593                     & \multicolumn{1}{c|}{0.0845} & \multicolumn{1}{l|}{0.0056131}             & \multicolumn{1}{l|}{0.0405}                    & 0.2482                                    \\
			
			\multicolumn{1}{c|}{1.5(0.5)}     & 0.0465                     & 0.1105                     & \multicolumn{1}{c|}{0.1352}  & 0.0471  & 0.0640                     & \multicolumn{1}{c|}{0.0909} & \multicolumn{1}{l|}{0.0036438}             & \multicolumn{1}{l|}{0.0395}                    & 0.2501                                    \\
			
			\multicolumn{1}{c|}{1.333(0.333)}   & 0.0400                     & 0.0679                    & \multicolumn{1}{c|}{0.1064}  & 0.0292 & 0.0350                    & \multicolumn{1}{c|}{0.0627} & \multicolumn{1}{l|}{0.0033354}             & \multicolumn{1}{l|}{0.0397}                    & 0.2432                                    \\
			
			\multicolumn{1}{c|}{1.25(0.25)}    & 0.0413                     & 0.0754                     & \multicolumn{1}{c|}{0.0878}  & 0.0229 & 0.0445                     & \multicolumn{1}{c|}{0.0580} & \multicolumn{1}{l|}{0.0009903}             & \multicolumn{1}{l|}{0.0436}                    & 0.2447                                    \\
			
			\multicolumn{1}{c|}{1.2(0.2)}     & 0.0388                     & 0.0716                     & \multicolumn{1}{c|}{0.1090} & 0.0196 & 0.0301                     & \multicolumn{1}{c|}{0.0456} & \multicolumn{1}{l|}{0.0007893}             & \multicolumn{1}{l|}{0.0358}                    & 0.2468                                    \\
			
			\multicolumn{1}{c|}{1.182(0.182)} & 0.0335                     & 0.0620                    & \multicolumn{1}{c|}{0.0894}  & 0.0216 & 0.0284                    & \multicolumn{1}{c|}{0.0365} & \multicolumn{1}{l|}{0.0007318}             & \multicolumn{1}{l|}{0.0335}                    & 0.2453                                    \\
			
			\multicolumn{1}{c|}{1.167(0.167)}   & 0.0304                     & 0.0708                     & \multicolumn{1}{c|}{0.1008}  & 0.0160 & 0.0341                     & \multicolumn{1}{c|}{0.0476} & \multicolumn{1}{l|}{0.0005918}             & \multicolumn{1}{l|}{0.0372}                    & 0.2462                                    \\
			
			\multicolumn{1}{c|}{1.1(0.1)}    & 0.0356                     & 0.0788                     & \multicolumn{1}{c|}{0.1001}  & 0.0191 & 0.0384                     & \multicolumn{1}{c|}{0.0489} & \multicolumn{1}{l|}{0.00048686}            & \multicolumn{1}{l|}{0.0332}                    & 0.2454                                    \\
			
			\multicolumn{1}{c|}{1.067(0.067)}  & 0.0358                     & 0.0652                     & \multicolumn{1}{c|}{0.1013}  & 0.0169 & 0.0318                     & \multicolumn{1}{c|}{0.0455} & \multicolumn{1}{l|}{0.00024720}            & \multicolumn{1}{l|}{0.0313}                    & 0.2445                                    \\
			
			\multicolumn{1}{l|}{1.05(0.05)}   & \multicolumn{1}{l}{0.0308} & \multicolumn{1}{l}{0.0625} & \multicolumn{1}{l|}{0.0845}  & 0.0149 & \multicolumn{1}{l}{0.0238} & \multicolumn{1}{l|}{0.0315} & \multicolumn{1}{l|}{0.00022473}            & \multicolumn{1}{l|}{0.0351}                    & 0.2443                                    \\ \hline
	\end{tabular}}
\end{table}
 \vfill\eject
 
 \bigskip\noindent
 {\bf Acknowledgement.} This research was supported by the Ministry of Human Capacities, Hungary grant TUDFO/47138-1/2019-ITM.

\end{document}